\documentclass[a4paper,11pt,reqno]{amsart}
\usepackage{amsmath}
\usepackage{amstext}
\usepackage{amsbsy}
\usepackage{amsopn}
\usepackage{upref}
\usepackage{amsthm}
\usepackage{amsfonts}
\usepackage{amssymb}
\usepackage{mathrsfs}
\usepackage{times}
\allowdisplaybreaks
 \newtheorem{theorem}{Theorem}

 \newtheorem{lemma}[theorem]{Lemma}
\theoremstyle{definition}
 \newtheorem{definition}{Definition}
\theoremstyle{remark}
 
\newcommand{\ep}{\varepsilon}
\newcommand{\p}{\partial}

\newcommand{\re}{\operatorname{Re}}
\newcommand{\im}{\operatorname{Im}}
\begin{document}
\title[dispersive equation]{Third order semilinear dispersive equations \\ related to deep water waves}
\author[H.~Chihara]{Hiroyuki CHIHARA}
\address{Mathematical Institute,  
         Tohoku University, 
         Sendai 980-8578, Japan}
\email{chihara@math.tohoku.ac.jp}
\begin{abstract}
We present local existence theorem of 
the initial value problem for 
third order semilinear dispersive partial differential equations in two space dimensions. 
This type of equations arises in the study of gravity wave of deep water, 
and cannot be solved by the classical energy method.  
To solve the initial value problem, 
we make full use of pseudodifferential operators 
with nonsmooth coefficients. 
\end{abstract}
\maketitle
\section{Introduction}
\label{section:1}
In this paper we study the initial value problem for 
third order semilinear dispersive equations of the form 
\begin{alignat}{2}
\p_tu+p(\p)u&=\sum_{j=0}^3a_jf_j(u)
&
\quad\text{in}\quad
&
\mathbb{R}^{1+2},
\label{equation:pde}
\\
u(0,x)&=u_0(x)
&
\quad\text{in}\quad
&
\mathbb{R}^2,
\label{equation:data}
\end{alignat}
where 
$u(t,x)$ is a complex-valued unknown function of 
$(t,x)=(t,x_1,x_2)\in\mathbb{R}^{1+2}$, 
$u_0(x)$ is an initial data, 
$\p_t=\p/\p{t}$, 
$\p_j=\p/\p{x_j}$, 
$\p=(\p_1,\p_2)$, 
$p(i\xi)$ is a pure-imaginary polynomial of degree three 
of $\xi=(\xi_1,\xi_2)\in\mathbb{R}^2$, 
$i=\sqrt{-1}$, 
$a_j\in\mathbb{C}$, 
$$
f_0(u)=uR_1\p_1\lvert{u}\rvert^2, 
f_1(u)=\lvert{u}\rvert^2\p_1u, 
f_2(u)=u^2\p_1\bar{u}, 
f_3{u}=\lvert{u}\rvert^2u, 
$$
$R_1=\p_1(-\Delta)^{-1/2}$ 
and 
$\Delta=\p_1^2+\p_2^2$. 
\par
This type of partial differential equations arises in 
the study of gravity wave of deep water. 
In fact, 
\begin{align*}
& \left(
  \p_t
  -
  \frac{1}{16}(\p_1^3-6\p_1\p_2^2)
  +
  \frac{i}{8}(\p_1^2-2\p_2^2)
  +
  \frac{1}{2}\p_1
  \right)u
\\
& \qquad\qquad
  =
  -
  \frac{i}{2}f_0(u)
  -
  \frac{3}{2}f_1(u)
  +
  \frac{1}{4}f_2(u)
  -
  \frac{i}{2}f_3(u)
\end{align*}
was derived by Dysthe in \cite{dysthe}, and 
\begin{align*}
& \left(
  \p_t
  -
  (b_1\p_1^3+b_2\p_1\p_2^2)
  +
  i(b_3\p_1^2+b_4\p_2^2)
  +
  b_5\p_1
  \right)u
\\
& \qquad\qquad
  =
  -
  \frac{i}{2}f_0(u)
  +
  \mu_1f_1(u)
  +
  \mu_2f_2(u)
  +
  i\mu_3f_3(u),
\end{align*}
$$
b_1, 
b_2, 
b_3, 
b_4, 
b_5, 
\mu_1, 
\mu_2, 
\mu_3 
\in 
\mathbb{R},   
$$
was formulated by Hogan in \cite{hogan}. 
Note that $a_0$ and $a_3$ are pure-imaginary, 
and $a_1$ and $a_2$ are real in 
the above physical models.
If $a_0\ne0$ or $\im{a_1}\ne0$, 
then the loss of derivatives occurs 
in \eqref{equation:pde}, 
and the classical energy method does not work. 
More precisely, 
$(\re{a_0})\lvert{u}\rvert^2R_1\p_1u$ in $a_0f_0(u)$,  
$a_0u^2R_1\p_1\bar{u}$ in $a_0f_0(u)$ 
and 
$(\im{a_1})\lvert{u}\rvert^2\p_1u$ in $a_1f_1(u)$ 
cannot be controled by the classical energy method. 
$f_2(u)$ has no problem 
since $f_2(u)\bar{u}=u^2\p_1(\bar{u}^2)/2$, 
which is controled by the integration by parts. 
\par
For this reason, there are few results on the existence of solutions to 
\eqref{equation:pde}-\eqref{equation:data}. 
In \cite{debouard}, 
using the abstract Cauchy-Kowalewski theorem, 
de Bouard proved time-local existence of a unique solution 
to the generalized equations of \eqref{equation:pde} 
for a real-analytic initial data. 
Recently, 
the time decay of the fundamental solution $e^{-tp(\p)}$ 
was discussed by Ben-Artzi, Koch and Saut in \cite{bks}. 
\par
Let $p_0(\p)$ be the principal part of $p(\p)$, 
that is, 
$p_0(\xi)$ is a homogeneous polynomial of degree three, 
and $p(\xi)-p_0(\xi)$ is a polynomial of degree two. 
It is said that 
$e^{-tp(\p)}$  has local smoothing effect 
if $e^{-tp(\p)}$ gains one derivative in $x$ 
locally in $\mathbb{R}^{1+2}$. 
More generelly, 
in case that $p(\p)$ is an operator of order $m>1$ on $\mathbb{R}^n$, 
it is said that $e^{-tp(\p)}$ has local smoothing effect 
if $e^{-tp(\p)}$ gains $(-\Delta)^{(m-1)/4}$ 
locally in $\mathbb{R}^{1+n}$. 
In \cite{hoshiro}, 
Hoshiro proved that 
$e^{-tp(\p)}$ has local smoothing effect if and only if 
\begin{equation}
p_0^\prime(\xi)=\nabla_\xi{p_0(\xi)}\ne0
\quad\text{for}\quad 
\xi\ne0.
\label{equation:nontrapping}
\end{equation}
\eqref{equation:nontrapping} is equivalent to so-called nontrapping condition of classical orbits
$$
\lim_{t\rightarrow\pm\infty}
\lvert{x+tp_0^\prime(\xi)}\rvert
=+\infty 
\quad\text{for}\quad
(x,\xi)\in\mathbb{R}\times\mathbb{R}\setminus\{0\}.
$$
For the detail of the smoothing effect of dispersive equations, see 
\cite{chihara2}, 
\cite{doi2}, 
\cite{hoshiro}, 
\cite{rs} 
and references therein. 
Under the nontrapping condition, 
one can expect the standard existence theorem for 
\eqref{equation:pde}-\eqref{equation:data}. 
\par
The aim of this paper is to present 
time-local existence theorem of 
\eqref{equation:pde}-\eqref{equation:data} 
in appropriate Sobolev spaces. 
Here we introduce notation regarding function spaces. 
Let $L^2(\mathbb{R}^2)$ be the set of 
all square-integrable functions on $\mathbb{R}^2$. 
For $f,g{\in}L^2(\mathbb{R}^2)$, set
$$
(f,g)
=
\int_{\mathbb{R}^2}f(x)\overline{g(x)}dx, 
\quad
\lVert{f}\rVert=\sqrt{(f,g)}. 
$$
Let $\langle{D}\rangle=(1-\Delta)^{1/2}$. 
For $s\in\mathbb{R}$, set 
$H^s(\mathbb{R}^2)=\langle{D}\rangle^{-s}L^2(\mathbb{R})$ 
and 
$\lVert{f}\rVert_s=\lVert\langle{D}\rangle^s{f}\rVert$. 
Let $I$ be an interval in $\mathbb{R}$. 
$C(I;H^s(\mathbb{R}^2))$ is the set of all 
$H^s(\mathbb{R}^2)$-valued continuous functions on $I$. 
$L^1(I;H^s(\mathbb{R}^2))$ is the set of all 
$H^s(\mathbb{R}^2)$-valued integrable functions on $I$. 
$L^\infty(I;H^s(\mathbb{R}^2))$ is the set of all 
$H^s(\mathbb{R}^2)$-valued 
essentially bounded functions on $I$. 
In this paper 
$(\cdot,\cdot)$ 
and 
$\lVert\cdot\rVert$ 
sometimes mean the inner product and the norm of 
$\mathbb{C}^l$-valued functions respectively, that is, 
$$
(U,V)=\sum_{\nu=1}^l(u_\nu,v_\nu), 
\quad
\lVert{U}\rVert^2
=
(U,U)
$$
for 
$
U=(u_1,\dotsc,u_l),
V=(v_1,\dotsc,v_l)
{\in}(L^2(\mathbb{R}^2))^l
$. 
Any confusion will not occur. 
Here we state our results. 
\begin{theorem}
\label{theorem:thm}
Suppose {\rm \eqref{equation:nontrapping}} and $s>3$. 
Then, for any $u_0{\in}H^s(\mathbb{R}^2)$, 
there exists $T=T(\lVert{u_0}\rVert_s)>0$ such that 
{\rm \eqref{equation:pde}-\eqref{equation:data}} 
possesses a unique solution 
$u{\in}C([-T,T];H^s(\mathbb{R}^2))$.  
\end{theorem}
Our method of proof of Theorem~\ref{theorem:thm} 
is an energy method via pseudodifferential calculus. 
This method was developed in \cite{chihara1} 
for semilinear Schr\"odinger equations of the form
\begin{equation}
\p_tu-i\Delta{u}=f(u,\p{u}) 
\quad\text{in}\quad
\mathbb{R}^{1+n}.
\label{equation:sch}
\end{equation}
The basic idea comes from the theory of well-posedness of 
the initial value problem for linear dispersive equations. 
See, e.g., \cite{chihara2}, \cite{doi1}, \cite{tarama1}, 
\cite{tarama2} and references therein. 
In particular, 
Tarama obtained the necessary and sufficient condition 
for the $L^2$-well-posedness of 
the initial value problem for 
one dimensional third order equations. 
See \cite{tarama1} and \cite{tarama2} for the detail. 
To apply linear theory to nonlinear equations including 
$\p{u}$, $\p\bar{u}$ and so on, 
it is very convenient to consider the system of 
the original equation and its complex conjugate. 
In other words, we consider the system for 
${}^t[u,\bar{u}]$, 
where $u$ is an unknown function of the original nonlinear equation. 
To control the bad first order terms 
by the local smoothing effect, 
we make use of a pseudodifferential operator of order zero 
discovered by Doi in \cite{doi1}. 
Using the unknown function $u$, 
we construct a nonsmooth symbol of 
the pseudodifferential operator. 
It is very important that 
the symbol is the type $(\rho,\delta)=(1,0)$, and 
the corresponding operator is easy to handle.  
\eqref{equation:pde}-\eqref{equation:data} is easier 
than the initial value problem for \eqref{equation:sch} 
since the local smoothing effect of $e^{-tp(\p)}$ 
is stronger than that of $e^{it\Delta}$, 
and the principal part of \eqref{equation:pde} 
and that of its complex conjugate are exactly same. 
In case of the semilinear Schr\"odinger equation 
\eqref{equation:sch}, 
the principal parts of equations of $u$ and $\bar{u}$ 
are $-i\Delta$ and $i\Delta$ respectively, 
and the diagonalization technique was used so that 
the smoothing estimates via pseudodifferential calculus 
worked. 
See \cite{chihara1} for the detail.  
If the principal part of \eqref{equation:sch} 
is not elliptic, then the diagonalization does not work. 
In this case, an operator of the exotic class was used, 
and higher smoothness of the initial data was required. 
See \cite{kpv} for the detail. 
Our symbolic calculus is based on 
the symbol smoothing method for 
pseudodifferential operators with nonsmooth coefficients. 
This method was discovered by Nagase in \cite{nagase}. 
For the basic pseudodifferential calculus, 
see e.g., 
\cite{kumanogo}, 
\cite{taylor1}, 
\cite{taylor2} 
and references therein. 
\par
The organization of this paper is as follows. 
In Section~\ref{section:2} 
we construct a sequence of approximate solutions 
by parabolic regularization. 
In Section~\ref{section:3} 
we summarize pseudodifferential calculus and related nonlinear estimates needed later. 
In Section~\ref{section:4} 
we study a linear system related to \eqref{equation:pde}. 
In Section~\ref{section:5} 
we prove Theorem~\ref{theorem:thm}.  
\section{A sequence of approximate solutions}
\label{section:2}
In this section 
we construct a sequence of approximate solutions 
by parabolic regularization. 
Let $\ep$ be a positive parameter. 
Consider a parametrized initial value problem of the form 
\begin{alignat}{2}
\p_tu+p(\p)u-\ep\Delta{u}&=\sum_{j=0}^3a_jf(u)
&
\quad\text{in}\quad
&
(0,+\infty)\times\mathbb{R}^2,
\label{equation:pdee}
\\
u(0,x)&=u_0(x)
&
\quad\text{in}\quad
&
\mathbb{R}^2. 
\label{equation:datae}
\end{alignat}
Since $e^{t(-p(\p)+\ep\Delta)}$ can control 
first order terms locally in time, 
\eqref{equation:pdee}-\eqref{equation:datae} 
has a time-local unique solution. 
\begin{lemma}
\label{theorem:parabolic} 
Suppose $s>2$. 
Then, for any $u_0{\in}H^s(\mathbb{R}^2)$, 
there exists $T_\ep=T(\ep,\lVert{u_0}\rVert_s)>0$ 
such that 
{\rm \eqref{equation:pdee}-\eqref{equation:datae}} 
possesses a unique solution 
$u{\in}C([0,T_\ep];H^s(\mathbb{R}^2))$. 
Moreover, the map $u_0\mapsto{u}$ is continuous.   
\end{lemma}
To prove Lemma~\ref{theorem:parabolic}, 
we need the following estimates. 
\begin{lemma}
\label{theorem:composition1}
{\rm (i)}\quad 
Let $s\in\mathbb{R}$. 
Then, for any $u{\in}H^s(\mathbb{R}^2)$ and $t>0$, 
$$
\lVert{e^{t(-p(\p)+\ep\Delta)}u}\rVert_s
\leqslant
C\left(1+\frac{1}{\sqrt{t\ep}}\right)
\lVert{u}\rVert_{s-1}.
$$
{\rm (ii)}\quad
Let $s>2$. 
Then, 
for any $j=0,1,2,3$ and $u,v\in{H^s(\mathbb{R}^2)}$,  
\begin{align*}
  \lVert{f_j(u)}\rVert_{s-1}
& \leqslant
  C\lVert{u}\rVert_s^3,
\\
  \lVert{f_j(u)-f_j(v)}\rVert_{s-1}
& \leqslant
  C
  (\lVert{u}\rVert_s^2+\lVert{v}\rVert_s^2)
  \lVert{u-v}\rVert_s.
\end{align*}
\end{lemma}
Here we introduce notation. 
For a multi-index $\alpha=(\alpha_1,\alpha_2)$, 
set $\lvert\alpha\rvert=\alpha_1+\alpha_2$, 
$\p^\alpha=\p_1^{\alpha_1}\p_2^{\alpha_2}$ 
and 
$\xi^\alpha=\xi_1^{\alpha_1}\xi_2^{\alpha_2}$. 
Let $\sigma\geqslant0$, 
and let $[\sigma]$ be the largest integer 
less than or equal to $\sigma$. 
$\mathscr{B}^\sigma(\mathbb{R}^2)$ is the set of all 
$C^{[\sigma]}$-functions $f(x)$ on $\mathbb{R}^2$ 
satisfying 
$\lVert{f}\rVert_{\mathscr{B}^\sigma}<+\infty$, where 
$$
\lVert{f}\rVert_{\mathscr{B}^\sigma}
=
\begin{cases}
\displaystyle\sup_{x\in\mathbb{R}^2}
\displaystyle\sum_{\lvert\alpha\rvert\leqslant{\sigma}}
\lvert\p^\alpha{f(x)}\rvert
& 
(\sigma=0,1,2,3,\dotsc)
\\
\lVert{f}\rVert_{\mathscr{B}^{[\sigma]}}
+
\displaystyle\sup_{\substack{x,y\in\mathbb{R}^2 \\ x\ne{y}}}
\displaystyle\sum_{\lvert\alpha\rvert=[\sigma]}
\dfrac{\lvert\p^\alpha{f(x)}-\p^\alpha{f(y)}\rvert}
     {\lvert{x-y}\rvert^{s-[s]}}
&
(\text{otherwise}) 
\end{cases}
$$
The Fourier transform of $u(x)$ is defined by 
$$
\hat{u}(\xi)
=
\frac{1}{2\pi}
\int_{\mathbb{R}^2}
e^{-ix\cdot\xi}u(x)dx,
\quad
x\cdot\xi=x_1\xi_1+x_2\xi_2.
$$
\begin{proof}[Proof of Lemma~\ref{theorem:composition1}]
(i)\quad
Set 
$\lvert\xi\rvert=\sqrt{\xi_1^2+\xi_2^2}$ 
and 
$\langle\xi\rangle=\sqrt{1+\lvert\xi\rvert^2}$ 
for short. 
In view of the Plancherel-Perseval formula, we deduce 
\begin{align*}
  \lVert{e^{t(-p(\p)+\ep\Delta)}u}\rVert_s
& =
  \lVert{\langle\xi\rangle^se^{t(-p(i\xi)-\ep\lvert\xi\rvert^2)}\hat{u}}\rVert
\\
& \leqslant
  \sup_{\xi\in\mathbb{R}^2}
  \langle\xi\rangle
  e^{-t\ep\lvert\xi\rvert^2}
  \lVert{u}\rVert_{s-1}
\\
& \leqslant
  \left(
        1
        +
        \sup_{\xi\in\mathbb{R}^2}
        \lvert\xi\rvert
        e^{-t\ep\lvert\xi\rvert^2}
   \right)
   \lVert{u}\rVert_{s-1}
\\
& \leqslant
  C\left(
        1+
        \frac{1}{\sqrt{t\ep}}
   \right)
   \lVert{u}\rVert_{s-1}.
\end{align*}
(ii)\quad
We show only the estimate of $f_0(u)$. 
Recall the Sobolev embedding 
$H^s(\mathbb{R}^2)\subset\mathscr{B}^\sigma(\mathbb{R}^2)$
for $s>1+\sigma$. 
It is easy to see that for $s>2$ and $u,v{\in}H^{s-1}(\mathbb{R})$, 
$$
\lVert{uv}\rVert_{s-1}
\leqslant
2^{s-1}
(\lVert{u}\rVert_{s-1}\lVert{v}\rVert_{\mathscr{B}^0}
 +
 \lVert{u}\rVert_{\mathscr{B}^0}\lVert{v}\rVert_{s-1})
\leqslant
C_s\lVert{u}\rVert_{s-1}\lVert{v}\rVert_{s-1}.
$$
Using this estimate and the $L^2$-boundedness of $R_1$, 
we deduce 
\begin{align*}
  \lVert{f_0(u)}\rVert_{s-1}
& \leqslant
  C_s
  \lVert{u}\rVert_{s-1}
  \rVert{R_1\p_1\lvert{u}\rvert^2}\rVert_{s-1}
\\
& \leqslant
  C_s
  \lVert{u}\rVert_{s-1}
  \rVert{\p_1\lvert{u}\rvert^2}\rVert_{s-1}
\\
& \leqslant
  2C_s
  \lVert{u}\rVert_{s-1}
  \rVert{\bar{u}\p_1u}\rVert_{s-1}
\\
& \leqslant
  2C_s^2
  \lVert{u}\rVert_{s-1}^2
  \lVert{\p_1u}\rVert_{s-1}
\\
& \leqslant
  2C_s^2
  \lVert{u}\rVert_{s-1}^2
  \lVert{u}\rVert_s.
\end{align*}
The other estimates can be obtained 
by similar computation. 
We omit the detail.
\end{proof} 
Applying the contraction mapping theorem and 
Lemma~\ref{theorem:composition1} to 
an integral equation
\begin{equation}
u(t)
=
e^{t(-p(\p)+\ep\Delta)}u_0
+
\int_0^t
e^{(t-\tau)(-p(\p)+\ep\Delta)}
\sum_{j=0}^3
a_jf_j(u(\tau))
d\tau, 
\label{equation:integraleqn}
\end{equation}
we can show Lemma~\ref{theorem:parabolic}. 
We omit the detail. 
\section{Pseudodifferential operators with nonsmooth coefficients}
\label{section:3}
Following \cite{chihara1}, 
we give the pseudodifferential calculus used later. 
We introduce symbols with nonsmooth coefficients. 
\begin{definition}
\label{definition:symbol}
Let $m\in\mathbb{R}$ and $\sigma\geqslant0$. 
$\mathscr{B}^\sigma\Psi^m(\mathbb{R}^n)$ 
is the set of all functions $p(x,\xi)$ on 
$\mathbb{R}^n\times\mathbb{R}^n$ satisfying 
$$
\lVert{p}\rVert_{\mathscr{B}^\sigma\Psi^m,l}
=
\sup_{\xi\in\mathbb{R}^n}
\sum_{\lvert\alpha\rvert\leqslant{l}}
\langle\xi\rangle^{m-\lvert\alpha\rvert}
\lVert\p_\xi^\alpha{p(\cdot,\xi)}\rVert_{\mathscr{B}^\sigma}
<+\infty
$$
for $l=0,1,2,3,\dotsc$.
\end{definition}
Set $D=-i\p$. 
For $p(x,\xi)\in\mathscr{B}^\sigma\Psi^m(\mathbb{R}^n)$, 
a pseudodifferential operator $p(x,D)$ is defined by 
$$
p(x,D)u(x)
=
\frac{1}{(2\pi)^n}
\iint_{\mathbb{R}^n\times\mathbb{R}^n}
e^{i(x-y)\cdot\xi}p(x,\xi)u(y)dyd\xi,
$$
where 
$x\cdot\xi=x_1\xi_1+\dotsb+x_n\xi_n$ 
for 
$x=(x_1,\dotsb,x_n)$ 
and 
$\xi=(\xi_1,\dotsc,\xi_n)$. 
Conversely, if an operator $P$ is given, 
then, its symbol is given by 
$\sigma(P)(x,\xi)=e^{-ix\cdot\xi}Pe^{ix\cdot\xi}$.  
The properties of the above 
pseudodifferential operators are the following. 
\begin{lemma}
\label{theorem:nagase}
{\rm (i)}\quad
Suppose $\sigma>0$ and 
$p(x,\xi)\in\mathscr{B}^\sigma\Psi^0(\mathbb{R}^n)$. 
Then, there exist 
$l\in\mathbb{N}$ and $C>0$ such that 
for any $u{\in}L^2(\mathbb{R}^n)$, 
$$
\lVert{p(x,D)u}\rVert
\leqslant
C
\lVert{p}\rVert_{\mathscr{B}^\sigma\Psi^0,l}
\lVert{u}\rVert.
$$
{\rm (ii)}\quad
Suppose $\sigma>1$ and 
$p_m(x,\xi)\in\mathscr{B}^\sigma\Psi^m(\mathbb{R}^n)$, 
{\rm ($m=0,1$)}. 
Set 
$$
q(x,\xi)=p_0(x,\xi)p_1(x,\xi), 
\quad
r(x,\xi)=\overline{p_1(x,\xi)}.
$$
Let $p_1(x,D)^\star$ be the formal adjoint of $p(x,D)$. 
Then, there exist 
$l\in\mathbb{N}$ and $C>0$ such that 
for any $u{\in}L^2(\mathbb{R}^n)$, 
\begin{align*}
  \lVert{(p_0(x,D)p_1(x,D)-q(x,D))u}\rVert
& \leqslant
  C
  \lVert{p_0}\rVert_{\mathscr{B}^\sigma\Psi^0,l}
  \lVert{p_1}\rVert_{\mathscr{B}^\sigma\Psi^1,l}
  \lVert{u}\rVert,
\\
  \lVert{(p_1(x,D)p_0(x,D)-q(x,D))u}\rVert
& \leqslant
  C
  \lVert{p_0}\rVert_{\mathscr{B}^\sigma\Psi^0,l}
  \lVert{p_1}\rVert_{\mathscr{B}^\sigma\Psi^1,l}
  \lVert{u}\rVert,
\\
  \lVert{(p_1(x,D)^\star-r(x,D))u}\rVert
& \leqslant
  C
  \lVert{p_1}\rVert_{\mathscr{B}^\sigma\Psi^1,l}
  \lVert{u}\rVert. 
\end{align*}
{\rm (iii)}\quad
Let $p(x,\xi)=[p_{jk}(x,\xi)]_{j,k=1,\dotsc,l}$ 
be an $l{\times}l$ matrix. 
Suppose that  
$p_{jk}(x,\xi)\in\mathscr{B}^2\Psi^1(\mathbb{R}^n)$, 
and 
there exists $\lambda\geqslant0$ such that 
$$
p(x,\xi)+{}^t\overline{p(x,\xi)}\geqslant0
\quad\text{for}\quad
\lvert\xi\rvert\geqslant{\lambda}.
$$
Then, there exist 
$l\in\mathbb{N}$ and $C>0$ such that 
for any $U{\in}(L^2(\mathbb{R}^n))^l$, 
\begin{equation}
\re(p(x,D)U,U)
\geqslant
-C
\sum_{j,k=1}^l
\lVert{p_{jk}}\rVert_{\mathscr{B}^2\Psi^1,l}
\lVert{U}\rVert^2.
\label{equation:garding}
\end{equation}
\end{lemma}
For the proof of Lemma~\ref{theorem:nagase}, 
see \cite{chihara1} and references therein. 
\par
In the proof of Theorem~\ref{theorem:thm}, 
we construct a symbol by using functions of the form 
\begin{align*}
  \phi_{1,\ep}(t,x_1)
& =
  \int_{\mathbb{R}}
  \lvert
  \langle{D_2}\rangle^{1/2+\sigma}
  u_{\ep}(t,x_1,x_2)
  \rvert^2
  dx_2,
\\
  \phi_{2,\ep}(t,x_2)
& =
  \int_{\mathbb{R}}
  \lvert
  \langle{D_1}\rangle^{1/2+\sigma}
  u_{\ep}(t,x_1,x_2)
  \rvert^2
  dx_1,
\end{align*}
where $\sigma$ is a small positive number, 
$\langle{D_j}\rangle=(1-\p_j^2)^{1/2}$, 
and 
$u_\ep$ is a solution to 
\eqref{equation:pdee}-\eqref{equation:datae}. 
We will use the following properties of 
$\phi_{1,\ep}$ and $\phi_{2,\ep}$. 
\begin{lemma}
\label{theorem:function} 
Let $s>3$, and 
let $u_\ep{\in}C([0,T_\ep];H^s(\mathbb{R}^2))$ 
be a solution to 
\eqref{equation:pdee}-\eqref{equation:datae}. 
Pick up $\sigma\in(0,1)$ so that $s\geqslant3+3\sigma$. 
Then, 
\begin{equation}
\phi_{1,\ep}(t,y), \phi_{2,\ep}(t,y) 
\in 
C([0,T_\ep];\mathscr{B}^{2+\sigma}(\mathbb{R})),
\label{equation:hoelder}
\end{equation}
\begin{equation}
\int_{-\infty}^y
\phi_{1,\ep}(t,z)
dz, 
\int_{-\infty}^y
\phi_{2,\ep}(t,z)
dz 
\in 
C^1([0,T_\ep];\mathscr{B}^\sigma(\mathbb{R})).
\label{equation:tbibun}
\end{equation}
\end{lemma}
\begin{proof}
Firstly, we evaluate 
$(1-\p_1^2)\phi_{1,\ep}(t,x_1)$. 
Using the Fourier inversion formula on $x_1$, 
the Plancherel-Perseval formula on $x_2$, 
and the Schwarz inequality on $\xi_1$, 
we deduce 
\begin{align}
& \lvert(1-\p_1^2)\phi_{1,\ep}(t,x_1)\rvert 
\nonumber
\\
& \leqslant
  2
  \int_{\mathbb{R}}
  \lvert
  \langle{D_1}\rangle^2
  \langle{D_2}\rangle^{1/2+\sigma}
  u_\ep(t,x_1,x_2)
  \rvert^2
  dx_2
\nonumber
\\
& =
  C
  \int_{\mathbb{R}}
  \left\lvert
  \int_{\mathbb{R}}
  e^{ix_1\xi_1}
  \langle\xi_1\rangle^2
  \langle\xi_2\rangle^{1/2+\sigma}
  \hat{u}_\ep(t,\xi_1,\xi_2)
  d\xi_1
  \right\rvert^2
  d\xi_2
\nonumber
\\
& \leqslant
  C
  \int_{\mathbb{R}}
  \left\lvert
  \int_{\mathbb{R}}
  \langle\xi_1\rangle^2
  \langle\xi_2\rangle^{1/2+\sigma}
  \lvert\hat{u}_\ep(t,\xi_1,\xi_2)\rvert
  d\xi_1
  \right\rvert^2
  d\xi_2
\nonumber
\\
& \leqslant
  C
  \int_{\mathbb{R}}
  \left(
  \int_{\mathbb{R}}
  \langle\xi_1\rangle^{-1-2\sigma}
  d\xi_1
  \right)
\nonumber
\\
& \qquad\qquad
  \times
  \left(
  \int_{\mathbb{R}}
  \lvert
  \langle\xi_1\rangle^{5/2+\sigma}
  \langle\xi_2\rangle^{1/2+\sigma}
  \hat{u}_\ep(t,\xi_1,\xi_2)
  \rvert^2
  d\xi_1
  \right)
  d\xi_2
\nonumber
\\
& \leqslant
  C
  \int_{\mathbb{R}}
  \int_{\mathbb{R}}
  \lvert
  \langle\xi_1\rangle^{5/2+\sigma}
  \langle\xi_2\rangle^{1/2+\sigma}
  \hat{u}_\ep(t,\xi_1,\xi_2)
  \rvert^2
  d\xi_1
  d\xi_2
\nonumber
\\
& \leqslant
  C
  \lVert{u_\ep}(t)\rVert_{3+2\sigma}^2.
\label{equation:301}
\end{align}
Secondly, we evaluate 
\begin{align*}
& \p_1^2\phi_{1,\ep}(t,x_1)-\p_1^2\phi_{1,\ep}(t,y_1)
\\
& =
  2\re
  \int_{\mathbb{R}}
  \langle\xi_2\rangle^{1/2+\sigma}
  \p_1^2u_\ep(t,x_1,x_2)
  \langle\xi_2\rangle^{1/2+\sigma}
  \overline{u_\ep(t,x_1,x_2)}
  dx_2
\\
& -
  2\re
  \int_{\mathbb{R}}
  \langle\xi_2\rangle^{1/2+\sigma}
  \p_1^2u_\ep(t,y_1,x_2)
  \langle\xi_2\rangle^{1/2+\sigma}
  \overline{u_\ep(t,y_1,x_2)}
  dx_2
\\
& +
  2
  \int_{\mathbb{R}}
  \langle\xi_2\rangle^{1/2+\sigma}
  \p_1u_\ep(t,x_1,x_2)
  \langle\xi_2\rangle^{1/2+\sigma}
  \p_1\overline{u_\ep(t,x_1,x_2)}
  dx_2
\\
& -
  2
  \int_{\mathbb{R}}
  \langle\xi_2\rangle^{1/2+\sigma}
  \p_1u_\ep(t,y_1,x_2)
  \langle\xi_2\rangle^{1/2+\sigma}
  \p_1\overline{u_\ep(t,y_1,x_2)}
  dx_2.   
\end{align*}
Note that 
$
\lvert{e^{ix_1\xi_1}-e^{iy_1\xi_1}}\rvert
\leqslant
2
\lvert\xi_1\rvert^\sigma
\lvert{x_1-y_1}\rvert^\sigma
$. 
In the same way as \eqref{equation:301}, we deduce 
\begin{align}
& \lvert\p_1^2\phi_{1,\ep}(t,x_1)-\p_1^2\phi_{1,\ep}(t,y_1)\rvert 
\nonumber
\\
& \leqslant
  C
  \max_{k=0,1,2}
  \int_{\mathbb{R}}
  \left\lvert
  \int_{\mathbb{R}}
  (e^{ix_1\xi_1}-e^{iy_1\xi_1})
  \xi_1^k
  \langle{\xi_2}\rangle^{1/2+\sigma}
  \hat{u}_\ep(t,\xi)
  d\xi_1
  \right\rvert
\nonumber
\\
& \qquad\qquad\times
  \left(
  \int_{\mathbb{R}}
  \lvert
  \langle{\xi}\rangle^{5/2+\sigma}
  \hat{u}_\ep(t,\xi)
  \rvert
  d\xi_1
  \right)
  d\xi_2
\nonumber
\\
& \leqslant
  C\lvert{x_1-y_1}\rvert^\sigma
  \int_{\mathbb{R}}
  \left(
  \int_{\mathbb{R}}
  \lvert
  \langle{\xi}\rangle^{5/2+2\sigma}
  \hat{u}_\ep(t,\xi)
  \rvert
  d\xi_1
  \right)^2
  d\xi_2
\nonumber
\\
& \leqslant
  C\lvert{x_1-y_1}\rvert^\sigma
  \int_{\mathbb{R}}
  \int_{\mathbb{R}}
  \lvert
  \langle{\xi}\rangle^{3+3\sigma}
  \hat{u}_\ep(t,\xi)
  \rvert^2
  d\xi_1
  d\xi_2
\nonumber
\\
& =
  C
  \lVert{u_\ep(t)}\rVert_{3+3\sigma}^2
  \lvert{x_1-y_1}\rvert^\sigma.  
\label{equation:302}
\end{align}
Combining \eqref{equation:301} and \eqref{equation:302}, 
we obtain \eqref{equation:hoelder}. 
\par
Next, we show \eqref{equation:tbibun}. 
Set $p_\ep(\p)=p(\p)-\ep\Delta$ for short. 
Note that $p_0(\xi)=p_0^\prime(\xi)\cdot\xi/3$.
Using the integration by parts, 
we deduce 
\begin{align*}
& \p_t
  \int_{-\infty}^{x_1}
  \phi_{1,\ep}(t,y_1)
  dy_1
\\
& =
  2\re
  \int_{-\infty}^{x_1}
  \int_{\mathbb{R}}
  \langle{D_2}\rangle^{1/2+\sigma}
  \p_tu_\ep(t,y_1,x_2)
  \langle{D_2}\rangle^{1/2+\sigma}
  \overline{u_\ep(t,y_1,x_2)}
  dy_1dx_2
\\
& =
  -2\re
  \int_{-\infty}^{x_1}
  \int_{\mathbb{R}}
  p_\ep(\p)
  \langle{D_2}\rangle^{1/2+\sigma}
  u_\ep(t,y_1,x_2)
  \langle{D_2}\rangle^{1/2+\sigma}
  \overline{u_\ep(t,y_1,x_2)}
  dy_1dx_2
\\
& +
  2\re
  \int_{-\infty}^{x_1}
  \int_{\mathbb{R}}
  \sum_{j=0}^3
  \langle{D_2}\rangle^{1/2+\sigma}
  a_jf_j(u_\ep)
  \langle{D_2}\rangle^{1/2+\sigma}
  \overline{u_\ep}
  dy_1dx_2
\\
& =
  -\frac{2}{3}\re
  \int_{-\infty}^{x_1}
  \int_{\mathbb{R}}
  p_0^\prime(\p)\cdot\p
  \langle{D_2}\rangle^{1/2+\sigma}
  u_\ep(t,y_1,x_2)
  \langle{D_2}\rangle^{1/2+\sigma}
  \overline{u_\ep(t,y_1,x_2)}
  dy_1dx_2
\\
& +
  2\re
  \int_{-\infty}^{x_1}
  \int_{\mathbb{R}}
  (p_\ep(\p)-p_0(\p))
  \langle{D_2}\rangle^{1/2+\sigma}
  u_\ep
  \langle{D_2}\rangle^{1/2+\sigma}
  \overline{u_\ep}
  dy_1dx_2
\\
& +
  2\re
  \int_{-\infty}^{x_1}
  \int_{\mathbb{R}}
  \sum_{j=0}^3
  \langle{D_2}\rangle^{1/2+\sigma}
  a_jf_j(u_\ep)
  \langle{D_2}\rangle^{1/2+\sigma}
  \overline{u_\ep}
  dy_1dx_2
\\
& =
  -\frac{2}{3}\re
  \int_{\mathbb{R}}
  \frac{\p p_0}{\p\xi_1}(\p)
  \langle{D_2}\rangle^{1/2+\sigma}
  u_\ep(t,x_1,x_2)
  \langle{D_2}\rangle^{1/2+\sigma}
  \overline{u_\ep(t,x_1,x_2)}
  dx_2
\\
& +
  \frac{2}{3}\re
  \int_{-\infty}^{x_1}
  \int_{\mathbb{R}}
  p_0^\prime(\p)
  \langle{D_2}\rangle^{1/2+\sigma}
  u_\ep
  \cdot
  \p
  \langle{D_2}\rangle^{1/2+\sigma}
  \overline{u_\ep}
  dy_1dx_2
\\
& +
  2\re
  \int_{-\infty}^{x_1}
  \int_{\mathbb{R}}
  (p_\ep(\p)-p_0(\p))
  \langle{D_2}\rangle^{1/2+\sigma}
  u_\ep
  \langle{D_2}\rangle^{1/2+\sigma}
  \overline{u_\ep}
  dy_1dx_2
\\
& +
  2\re
  \int_{-\infty}^{x_1}
  \int_{\mathbb{R}}
  \sum_{j=0}^3
  \langle{D_2}\rangle^{1/2+\sigma}
  a_jf_j(u_\ep)
  \langle{D_2}\rangle^{1/2+\sigma}
  \overline{u_\ep}
  dy_1dx_2.
\end{align*}
Then, we have 
\begin{align}
  \left\lvert
  \p_t
  \int_{-\infty}^{x_1}
  \phi_{1,\ep}(t,y_1)
  dy_1
  \right\rvert
& \leqslant
  C
  \int_{\mathbb{R}}
  \lvert
  \langle{D}\rangle^{5/2+\sigma}
  u_\ep(t,x_1,x_2)
  \rvert^2
  dx_2
\nonumber
\\
& +
  C
  \lVert{u_\ep(t)}\rVert_{5/2+\sigma}^2
  +
  C
  \sum_{j=0}^3
  \lVert{f_j(u_\ep(t))}\rVert_{1/2+\sigma}
  \lVert{u_\ep(t)}\rVert_{1/2+\sigma}
\nonumber
\\
& \leqslant
  C
  (
  \lVert{u_\ep(t)}\rVert_{3+2\sigma}^2
  +
  \lVert{u_\ep(t)}\rVert_{3+2\sigma}^4  
  ).
\label{equation:303}
\end{align}
Similarly, we can get
\begin{equation}
\left\lvert
\p_t
\int_{y_1}^{x_1}
\phi_{1,\ep}(t,z_1)
dz_1
\right\rvert
\leqslant
C
(
\lVert{u_\ep(t)}\rVert_{3+3\sigma}^2
+
\lVert{u_\ep(t)}\rVert_{3+3\sigma}^4  
)
\lvert{x_1-y_1}\rvert^\sigma.
\label{equation:304}
\end{equation}
Combining \eqref{equation:303} and \eqref{equation:304}, 
we obtain \eqref{equation:tbibun}.
\end{proof}
%
%
\section{A linear pseudodifferential system}
\label{section:4}
As is considered in \cite{chihara1}, 
it is very convenient to establish 
the $H^s$-well-posedness of the initial value problem 
for a $2\times2$ system related to the system for 
\eqref{equation:pde} and its complex conjugate. 
Consider the initial value problem of the form
\begin{alignat}{2}
LU&=F(t,x)
&
\quad\text{in}\quad
&
(0,T)\times\mathbb{R}^2,
\label{equation:system}
\\
U(0,x)&=U_0(x)
&
\quad\text{in}\quad
&
\mathbb{R}^2,
\label{equation:datas}
\end{alignat}
where 
$U(t,x)$ is a $\mathbb{C}^2$-valued unknown function, 
$F(t,x)$ and $U_0(x)$ are given 
$\mathbb{C}^2$-valued functions, 
$T>0$, 
$$
L
=
I(\p_t+p_0(\p)-\ep\Delta)
+iJp_1(\p)
+
A(t)
+
B(t),
$$
$p_0(\xi)$ is same as that of \eqref{equation:pde}, 
$p_1(\xi)$ is a homogeneous real polynomial of degree two, 
$\ep\geqslant0$,  
$$
I
=
\begin{bmatrix}
1 & 0
\\
0 & 1 
\end{bmatrix},
\quad
J
=
\begin{bmatrix}
1 & 0
\\
0 & -1
\end{bmatrix},
$$
$$
A(t)
=
\begin{bmatrix}
a_{11}(t,x,D) & a_{12}(t,x,D)
\\
a_{21}(t,x,D) & a_{22}(t,x,D)
\end{bmatrix},
B(t)
=
\begin{bmatrix}
B_{11}(t) & B_{12}(t)
\\
B_{21}(t) & B_{22}(t) 
\end{bmatrix},
$$
$
a_{jk}(t,x,\xi)
\in
C([0,T];\mathscr{B}^{2+\sigma}\Psi^1(\mathbb{R}^2))
$, 
($\sigma>0$), 
$B_{jk}(t)$ and 
$[\langle{D}\rangle^s,B_{jk}(t)]$, 
($s\in[0,1]$) are $L^2$-bounded. 
Here we state the definition of the well-posedness. 
\begin{definition}
\label{definition:hswp}
Let $s\in\mathbb{R}$. 
The initial value problem 
{\rm \eqref{equation:system}-\eqref{equation:datas}} 
is said to be $H^s$-well-posed if for any 
$U_0\in(H^s(\mathbb{R}^2))^2$ and 
$F\in(L^1(0,T;H^s(\mathbb{R}^2)))^2$, 
{\rm \eqref{equation:system}-\eqref{equation:datas}} 
possesses a unique solution 
$U\in(C([0,T];H^s(\mathbb{R}^2)))^2$. 
\end{definition}
We give a sufficient condition of $H^s$-well-posedness 
used later. 
\begin{lemma}
\label{theorem:syswp}
Suppose $s\in[0,1]$ and \eqref{equation:nontrapping}. 
If there exists 
$
\phi(t,y)
\in
C([0,T];\mathscr{B}^{2+\sigma}(\mathbb{R}))$ 
such that $\phi(t,y)\geqslant0$, 
$$
\sup_{t\in[0,T]}
\int_{\mathbb{R}}\phi(t,y)dy
+
\sup_{\substack{t\in[0,T] \\ x_1\in\mathbb{R}}}
\left\lvert
\int_{-\infty}^{x_1}
\p_t\phi(t,y)dy
\right\rvert
<+\infty,
$$
$$
\sum_{j=1,2}
\lvert\re{a_{jj}(t,x_1,x_2,\xi)}\rvert
+
\sum_{j+k=3}
\lvert{a_{jk}(t,x_1,x_2,\xi)}\rvert
\leqslant
\phi(t,x_\nu)\lvert\xi\rvert
$$
for $(t,x)=(t,x_1,x_2)\in[0,T]\times\mathbb{R}^2$, 
$\nu=1,2$ and $\lvert\xi\rvert\geqslant1$, 
then 
{\rm \eqref{equation:system}-\eqref{equation:datas}} 
is $H^2$-well-posed. 
\end{lemma}
\begin{proof}
Lemma~\ref{theorem:syswp} follows from energy estimates 
and duality argument. 
We show the forward $L^2$-estimate for $LU=F$, 
and omit the backward energy inequality 
for $L^\star{U}=F$.  
Set $K(t)=Ik(t,x,D)$, $k(t,x,\xi)=e^{\gamma(t,x,\xi)}$, 
$K^\prime(t)=Ik^\prime(t,x,,D)$, 
$k^\prime(t,x,\xi)=e^{-\gamma(t,x,\xi)}$, 
$$
\gamma(t,x,\xi)
=
\sum_{j=1,2}
\int_{-\infty}^{x_j}
\phi(t,y)dy
\frac{\p{p_0}}{\p\xi_j}(\xi)
\lvert{p_0^\prime(\xi)}\rvert^{-2}
\lvert\xi\rvert
\chi\left(\dfrac{\xi}{\lambda}\right), 
$$
where $\lambda$ is a positive constant determined later, 
and $\chi(\xi){\in}C^\infty(\mathbb{R})$ satisfying 
$\chi(\xi)=1$ for $\lvert\xi\rvert\geqslant1$ 
and 
$\chi(\xi)=0$ for $\lvert\xi\rvert\leqslant1/2$. 
Since 
$$
\sigma(K^\prime(t)K(t))(x,\xi), \ 
\sigma(K(t)K^\prime(t))(x,\xi)
=
I+O(\lambda^{-3}),
$$
there exists $\lambda_0>0$ such that 
$K(t)$ is invertible on $(L^2(\mathbb{R}^2))^2$ 
for $\lambda\geqslant\lambda_0$. 
Fix $\lambda\geqslant1+\lambda_0$. 
Note that 
$$
\sigma([k(t,x,D),p_0(\p)]-Q_0(t)k(t,x,D))(x,\xi)
\in
C([0,T];\mathscr{B}^{1+\sigma}\Psi^0(\mathbb{R}^2)),
$$
\begin{align*}
  \sigma(Q_0(t))(x,\xi)
& =
  \p{e^{\gamma(t,x,\xi)}}\cdot{p_0^\prime(\xi)}
  e^{-\gamma(t,x,\xi)}
\\
& =
  \sum_{j=1,2}
  \phi(t,x_j)
  \left\lvert\frac{\p{p_0}}{\p\xi_j}(\xi)\right\rvert^2
  \frac{\lvert\xi\rvert}{\lvert{p_0^\prime(\xi)}\rvert^2}
  \chi\left(\frac{\xi}{\lambda}\right).
\end{align*}
Set $Q(t)=IQ_0(t)+A(t)$ for short. 
Applying $K(t)$ to $L$, we deduce 
$$
K(t)L
=
\{
I(\p_t+p_0(\p)-\ep\Delta)+iJp_1(p)+Q(t)+R(t)
\}K(t),
$$
\begin{align}
  R(t)
& =
  \biggl\{
  \frac{\p{K}}{\p{t}}(t)
  +
  I\Bigl([k(t,x,D),p_0(\p)]-Q_0(t)k(t,x,D)\Bigr)
\nonumber
\\
& \qquad
  +
  iJ[k(t,x,D),p_1(\p)]
  -
  \ep{I}[k(t,x,D),\Delta]
\nonumber
\\
& \quad\qquad
  +
  \Bigl[[k(t,x,D),a_{jk}(t,x,D)]\Bigr]_{j,k=1,2}
  +
  K(t)B(t)
  \biggr\}
  K(t)^{-1}.
\label{equation:nanako} 
\end{align}
In view of the hypothesis, 
$$
\sigma(Q(t))(x,\xi)
+
\sigma(Q(t)^\star)(x,\xi)
\geqslant0
\quad\text{for}\quad
\lvert\xi\rvert\geqslant\lambda,
$$
and $R(t)$ is $L^2$-bounded. 
Then, we deduce 
\begin{align*}
  \frac{d}{dt}
  \lVert{K(t)U(t)}\rVert^2
& =
  2\re({\p_t}K(t)U(t),K(t)U(t))
\\
& =
  2\re((-Ip_0(\p)-iJp_1(\p))K(t)U(t),K(t)U(t))
\\
& +
  2\re((\ep{I}\Delta-Q(t))K(t)U(t),K(t)U(t))
\\
& -
  2\re(R(t)K(t)U(t),K(t)U(t))
  +
  2\re(K(t)F(t),K(t)U(t))
\\
& =
  -
  2\ep\lVert\p{K(t)U(t)}\rVert^2
  -
  2\re(Q(t)K(t)U(t),K(t)U(t))
\\
& -
  2\re(R(t)K(t)U(t),K(t)U(t))
  +
  2\re(K(t)F(t),K(t)U(t))
\\
& \leqslant
  -
  2\re(Q(t)K(t)U(t),K(t)U(t))
\\
& +
  2C\lVert{K(t)U(t)}\rVert^2
  +
  2\lVert{K(t)F(t)}\rVert\lVert{K(t)U(t)}\rVert. 
\end{align*}
By the sharp G{\aa}rding inequality 
\eqref{equation:garding}, 
$$
\frac{d}{dt}
\lVert{K(t)U(t)}\rVert^2
\leqslant
  2C_0\lVert{K(t)U(t)}\rVert^2
  +
  2\lVert{K(t)F(t)}\rVert\lVert{K(t)U(t)}\rVert,
$$
which becomes 
\begin{equation}
\frac{d}{dt}
\lVert{K(t)U(t)}\rVert
\leqslant
C_0\lVert{K(t)U(t)}\rVert
+
\lVert{K(t)F(t)}\rVert.
\label{equation:wnua}
\end{equation}
Integrating \eqref{equation:wnua} over $[0,t]$, we obtain 
$$
\lVert{K(t)U(t)}\rVert
\leqslant
e^{C_0t}\lVert{K(0)U_0}\rVert
+
\int_0^t
e^{C_0(t-\tau)}\lVert{K(\tau)F(\tau)}\rVert, 
$$
which is a desired energy inequality. 
\par
Let $s\in(0,1]$. 
Applying $I\langle{D}\rangle^s$ 
to \eqref{equation:system}, we have 
$$
(L+\tilde{A}_s(t))\langle{D}\rangle^sU
=
\tilde{B}_s(t)U+\langle{D}\rangle^sF(t,x),
$$
$$
\tilde{A}_s(t)
=
\begin{bmatrix}
[\langle{D}\rangle^s,a_{11}(t,x,D)]
&
[\langle{D}\rangle^s,a_{12}(t,x,D)]
\\
[\langle{D}\rangle^s,a_{21}(t,x,D)]
&
[\langle{D}\rangle^s,a_{22}(t,x,D)]
\end{bmatrix}
\langle{D}\rangle^{-s}, 
$$
$$
\tilde{B}_s(t)
=
-
\begin{bmatrix}
[\langle{D}\rangle^s,B_{11}(t)]
&
[\langle{D}\rangle^s,B_{12}(t)]
\\
[\langle{D}\rangle^s,B_{21}(t)]
&
[\langle{D}\rangle^s,B_{22}(t)]
\end{bmatrix}.
$$
$\tilde{A}_s(t)$ and $\tilde{B}_s(t)$ are $L^2$-bounded. 
If \eqref{equation:system}-\eqref{equation:datas} 
is $L^2$-well-posed, then 
The initial value problem for a modified system 
\begin{equation}
(L+\tilde{A}_s(t))U
=
F(t,x)
\label{equation:95.5}
\end{equation} 
is also $L^2$-well-posed.   
The $L^2$-well-posedness of 
\eqref{equation:system}-\eqref{equation:datas} and 
\eqref{equation:95.5}-\eqref{equation:datas}    
implies the $H^s$-well-posedness of 
\eqref{equation:system}-\eqref{equation:datas}. 
This completes the proof. 
\end{proof}
\section{Energy method}
\label{section:5}
In this section we prove Theorem~\ref{theorem:thm}. 
Suppose $s>3$. 
Let $u_\ep{\in}C([0,T_\ep];H^s(\mathbb{R}^2))$ be 
a unique solution to 
\eqref{equation:pdee}-\eqref{equation:datae}. 
Firstly, we construct a solution 
$u{\in}L^\infty(0,T;H^s(\mathbb{R}^2))$ with some $T>0$ 
by compactness argument. 
More precisely, we show that there exists $T>0$ such that 
$T_\ep\geqslant{T}$ for any $\ep>0$, and 
$\{u_\ep\}_{\ep>0}$ is bounded in 
$L^\infty(0,T;H^s(\mathbb{R}^2))$. 
\par
Let $\chi(\xi)$ be the same as in Section~\ref{section:4}. 
We split $R_1$ into the principal part and the remainder part by 
$$
R_1=r_0(D)+\tilde{R}_1, 
\quad
r_0(\xi)
\frac{i\xi}{\lvert\xi\rvert}\chi(\xi), 
\quad
\tilde{R}_1=R_1(1-\chi(D)).
$$
Then, 
$$
r_0(\xi)
\in
\Psi^0(\mathbb{R}^2)
=
\bigcap_{k=1}^\infty
\mathscr{B}^k\Psi^0(\mathbb{R}^2),
\quad
\lVert{\tilde{R}_1v}\rVert_{m}
\leqslant
C_m\lVert{v}\rVert
$$
for any $m\geqslant0$ and $v{\in}L^2(\mathbb{R}^2)$. 
Since $p(i\xi)$ is pure-imaginary, 
$$
p(\p)=p_0(\p)+ip_1(\p)+p_2(\p)+ip_3,
$$
where $p_j(\xi)$ is a homogeneous real polynomial of 
degree $3-j$. 
Set $\theta=s-[s]$ for short. 
Let $\alpha$ be a multi-index satisfying 
$\lvert\alpha\rvert=[s]$. 
Applying $\langle{D}\rangle^\theta\p^\alpha$ to 
\eqref{equation:pdee}, we have 
\begin{align*}
& (\p_t+p_0(\p)+ip_1(\p)-\ep\Delta)
  \langle{D}\rangle^\theta\p^\alpha{u_\ep}
\\
& \quad
  +
  a_{1,\ep}(t,x,D)
  \langle{D}\rangle^\theta\p^\alpha{u_\ep}
  +
  a_{2,\ep}(t,x,D)
  \overline{\langle{D}\rangle^\theta\p^\alpha{u_\ep}
}
  =f_{\ep,\alpha},
\end{align*}
\begin{align*}
  a_{1,\ep}(t,x,\xi)
& =
  -
  i\lvert{u_\ep(t,x)}\rvert^2
  (a_0r_0(\xi)+a_1)\xi_1
  +
  ip_2(\xi),
\\
  a_{2,\ep}(t,x,\xi)
& =
  -
  iu_\ep(t,x)^2
  (a_0r_0(\xi)+a_2)\xi_1,
\end{align*}
\begin{align*}
  f_{\ep,\alpha}
& =
  \biggl\{
  \sum_{j=0}^3
  a_j
  \langle{D}\rangle^\theta\p^\alpha{f_j(u_\ep)}
\\
& \qquad
  -
  \lvert{u_\ep}\rvert^2(a_0R_1+a_1)\p_1
  \langle{D}\rangle^\theta\p^\alpha{u_\ep}
  -
  u_\ep^2(a_0R_1+a_1)\p_1
  \langle{D}\rangle^\theta\p^\alpha{\bar{u}_\ep}
  \biggr\}
\\
& +
  a_0\lvert{u_\ep}\rvert^2\tilde{R}_1\p_1
  \langle{D}\rangle^\theta\p^\alpha{u_\ep}
  +
  a_0u_\ep^2\tilde{R}_1\p_1
  \langle{D}\rangle^\theta\p^\alpha{\bar{u}_\ep}
\\
& -
  [\langle{D}\rangle^\theta,\lvert{u_\ep}\rvert^2(a_0R_1+a_1)\p_1]\p^\alpha{u_\ep}
\\
& -
  [\langle{D}\rangle^\theta,u_\ep^2(a_0R_1+a_1)\p_1]\p^\alpha{\bar{u}_\ep}
  -
  ip_3\langle{D}\rangle^\theta\p^\alpha{u_\ep}.
\end{align*}
Since the highest order of differentiation in 
$f_{\ep,\alpha}$ is $s=\theta+\lvert\alpha\rvert$, 
$$
\lVert{f_{\ep,\alpha}(t)}\rVert
\leqslant
C
(
\lVert{u_\ep(t)}\rVert_s
+
\lVert{u_\ep(t)}\rVert_s^3
).
$$
If we set 
$$
U_{\ep,\alpha}
=
\begin{bmatrix}
\langle{D}\rangle^\theta\p^\alpha{u_\ep}
\\
\overline{\langle{D}\rangle^\theta\p^\alpha{u_\ep}}
\end{bmatrix},
\quad
U_{\alpha,0}
=
\begin{bmatrix}
\langle{D}\rangle^\theta\p^\alpha{u_0}
\\
\overline{\langle{D}\rangle^\theta\p^\alpha{u_\ep}}
\end{bmatrix},
\quad
F_{\ep,\alpha}
=
\begin{bmatrix}
f_{\ep,\alpha}
\\
\overline{f_{\ep,\alpha}}
\end{bmatrix},
$$
$$
\sigma(A_\ep(t))(x,\xi)
=
\begin{bmatrix}
a_{1,\ep}(t,x,\xi)
&
a_{2,\ep}(t,x,\xi)
\\
\overline{a_{2,\ep}(t,x,-\xi)}
&
\overline{a_{1,\ep}(t,x,-\xi)}
\end{bmatrix},
$$
then $U_{\ep,\alpha}$ solves
\begin{alignat}{2}
  \Bigl\{
  I(\p_t+p_0(\p)-\ep\Delta)+iJp_1(\p)+A_\ep(t)
  \Bigr\}U_{\ep,\alpha}
&=F_{\ep,\alpha}
&
\quad\text{in}\quad
&
(0,T_\ep)\times\mathbb{R}^2,
\label{equation:sysalpha}
\\
U_{\ep,\alpha}(0,x)&=U_{\alpha,0}(x)
&
\quad\text{in}\quad
&
\mathbb{R}^2.
\nonumber
\end{alignat}
\par
Let $\phi_{1,\ep}$ and $\phi_{2,\ep}$ be the functions 
introduced in Section~\ref{section:3}. 
Since 
$$
\lvert\re{a_{1,\ep}(t,x,\xi)}\rvert
+
\lvert{a_{2,\ep}(t,x,\xi)}\rvert
\leqslant
(2\lvert{a_0}\rvert+\lvert{a_1}\rvert+\lvert{a_2}\rvert)
\lvert{u_\ep(t,x)}\rvert^2\lvert\xi_1\rvert, 
$$
there exists $C_0>0$ which is independent of $\ep>0$, 
such that 
\begin{equation}
2\lvert\re{a_{1,\ep}(t,x,\xi)}\rvert
+
2\lvert{a_{2,\ep}(t,x,\xi)}\rvert
\leqslant
C_0
\lvert\xi\rvert
\min\{\phi_{1,\ep}(t,x_1),\phi_{2,\ep}(t,x_2)\}.
\label{equation:dominate}
\end{equation}
Set 
$$
\phi_\ep(t,y)
=
C_0
\phi_{1,\ep}(t,y)
+
C_0
\phi_{2,\ep}(t,y),
$$
$$
\gamma_{\ep}(t,x,\xi)
=
\sum_{j=1,2}
\int_{-\infty}^{x_j}\phi_{\ep}(t,y)dy
\frac{\p{p_0}}{\p\xi_j}(\xi)
\frac{\lvert\xi\rvert}{\lvert{p_0^\prime(\xi)}\rvert^2}
\chi(\xi),
$$
$$
\sigma(K_\ep(t))(x,\xi)
=
Ie^{\gamma_\ep(t,x,\xi)},
\quad
\sigma(K_\ep^\prime(t))(x,\xi)
=
Ie^{-\gamma_\ep(t,x,\xi)}.
$$
In view of 
\eqref{equation:hoelder}, 
\eqref{equation:tbibun} 
and 
\eqref{equation:dominate}, 
$\phi_\ep$ satisfies the conditions in 
Lemma~\ref{theorem:syswp} for \eqref{equation:sysalpha}. 
\par
We evaluate 
$$
N_\ep(t)
=
\sum_{\lvert\alpha\rvert=[s]}
\lVert{K_\ep(t)U_{\ep,\alpha}(t)}\rVert
+
\lVert{u_\ep(t)}\rVert_{s-1}.
$$
Since $N_\ep(0)$ is independent of $\ep>0$, 
set $M=N_\ep(0)$ for short. 
Here we introduce 
$$
T_\ep^\star
=
\sup
\{T>0 
\vert 
N_\ep(t)\leqslant2M 
\ \text{for}\ 
t\in[0,T]
\}.
$$
Since 
$$
\sigma(K_\ep^\prime(t)K_\ep(t))(x,\xi), 
\sigma(K_\ep(t)K_\ep^\prime(t))(x,\xi) 
= 
I+O(\langle\xi\rangle^{-3}), 
$$
there exists $C_M>1$ which is independent of $\ep>0$, 
such that 
$$
C_M^{-1}
\lVert{u_\ep(t)}\rVert_s
\leqslant
N_\ep(t)
\leqslant
C_M
\lVert{u_\ep(t)}\rVert_s
\quad\text{for}\quad
t\in[0.T_\ep^\star]. 
$$
Applying $K_\ep(t)$ to \eqref{equation:sysalpha}, 
we have 
$$
\Bigl\{
I(\p_t+p_0(\p)-\ep\Delta)
+iJ\p_1(\p)+Q_\ep(t)
\Bigr\}K_\ep(t)U_{\ep,\alpha}
+
R_\ep(t)U_{\ep,\alpha}
=
K_\ep(t)F_{\ep,\alpha},
$$
$$
Q_\ep(t)
=
IQ_{\ep,0}(t)+A_\ep(t),
$$
$$
\sigma(Q_{\ep,0}(t))(x,\xi)
=
\sum_{j=1,2}
\phi_{\ep}(t,x_j)
\left\lvert\frac{\p{p_0}}{\p\xi_j}(\xi)\right\rvert^2
\frac{\lvert\xi\rvert}{\lvert{p_0^\prime(\xi)}\rvert^2}
\chi(\xi),
$$
$R_\ep(t)$ corresponds to $R(t)K(t)$ 
in \eqref{equation:nanako}. 
In view of Lemmas~\ref{theorem:nagase} and 
\ref{theorem:function}, 
we get 
$$
\lVert{R_\ep(t)U_{\ep,\alpha}(t)}\rVert
\leqslant
CM^2N_\ep(t) 
\quad\text{for}\quad
t\in[0,T_\ep^\star]. 
$$
In the same way as the energy estimate in Section~\ref{section:4}, 
we deduce 
\begin{align*}
  \frac{d}{dt}
  \lVert{K_\ep(t)U_{\ep,\alpha}(t)}\rVert^2
& \leqslant
  -
  2\re
  (Q_\ep(t)K_\ep(t)U_{\ep,\alpha}(t),K_\ep(t)U_{\ep,\alpha}(t))
\\
& +
  2
  \Bigl(
  \lVert{R_\ep(t)U_{\ep,\alpha}(t)}\rVert
  +
  \lVert{K_\ep(t)F_{\ep,\alpha}(t)}\rVert
  \Bigr)
  \lVert{K_\ep(t)U_{\ep,\alpha}(t)}\rVert
\\
& \leqslant
  2C_1M^2N_\ep(t)^2 
\end{align*}
for $t\in[0,T_\ep^\star]$, 
where $C_1>0$ depends only on $M$. 
Then, we have 
\begin{equation}
\lVert{K_\ep(t)U_{\ep,\alpha}(t)}\rVert
\leqslant
\lVert{K_\ep(0)U_{\alpha,0}}\rVert
+
C_1M^2
\int_0^tN_\ep(\tau)d\tau. 
\label{equation:maki}
\end{equation}
Using \eqref{equation:integraleqn}, we get 
\begin{align}
  \lVert{u_\ep(t)}\rVert_{s-1}
& \leqslant
  \lVert{u_0}\rVert_{s-1}
  +
  \sum_{j=0}^3
  \lvert{a_j}\rvert
  \int_0^t
  \lVert{f_j(u(\tau))}\rVert_{s-1}
  d\tau
\nonumber
\\
& \leqslant
  \lVert{u_0}\rVert_{s-1}
  +
  C_2M^2
  \int_0^tN_\ep(\tau)d\tau.
\label{equation:maoko}
\end{align}
Combining \eqref{equation:maki} 
and \eqref{equation:maoko}, 
we obtain 
$$
N_\ep(t)
\leqslant
N_\ep(0)
+
C_3M^2
\int_0^tN_\ep(\tau)d\tau
$$
for $t\in[0,T_\ep^\star]$, 
where $C_3>0$ is independent of $\ep>0$. 
The Gronwall inequality implies that 
$$
N_\ep(t)
\leqslant
M\exp\Bigl(C_3M^2t\Bigr)
\quad\text{for}\quad
t\in[0,T_\ep^\star]. 
$$
If we set $t=T_\ep^\star$, then 
$2M\leqslant\exp(C_3M^2T_\ep^\star)$, 
which gives 
$T_\ep^\star\geqslant\log2/C_3M^2$. 
Set $T=\log2/C_3M^2$ for short. 
$\{u_\ep\}_{\ep>0}$ is bounded in 
$L^\infty(0,T;H^s(\mathbb{R}^2))$. 
The standard compactness argument shows that 
there exist a subsequence $\{u_\ep\}$ and $u$ such that 
\begin{alignat*}{2}
  u_\ep
& \longrightarrow
  u
&
  \quad\text{in}\quad
& L^\infty(0,T;H^s(\mathbb{R}^2)) 
  \quad\text{weakly}^\star,  
\\
  u_\ep
& \longrightarrow
  u
&
  \quad\text{in}\quad
& C([0,T];H^{s-\delta}_{\text{loc}}(\mathbb{R}^2)),
  \quad
  (\delta>0),  
\end{alignat*}
as $\ep\downarrow0$. 
It is easy to see 
\begin{equation}
u
\in
L^\infty(0,T;H^s(\mathbb{R}^2))
\cap
C([0,T];H^{s-\delta}(\mathbb{R}^2)), 
\quad
(\delta>0),
\label{equation:ryoko}
\end{equation}
and $u$ solves 
\eqref{equation:pde}-\eqref{equation:data} 
in the sense of distribution. 
\par
Secondly, we prove the uniqueness of solution. 
Let 
$u,v{\in}L^\infty(0,T;H^s(\mathbb{R}^2))$ 
be solutions to \eqref{equation:pde} 
with $u(0)=v(0)$. 
Set $w=u-v$ for short. 
Then, $w(0)=0$, and $w$ solves 
\begin{equation}
(\p_t+p(\p)+a_1(t,x,D)+B_1(t))w
+
(a_2(t,x,D)+B_2(t))\bar{w}=0, 
\label{equation:difference}
\end{equation}
\begin{align*}
  a_1(t,x,\xi)
& =
  -
  ia_0\lvert{v(t,x)}\rvert^2r_0(\xi)\xi_1
  -
  ia_1\lvert{v(t,x)}\rvert^2\xi_1,
\\
  a_2(t,x,\xi)
& =
  -
  ia_0u(t,x)v(t,x)r_0(\xi)\xi_1
  -
  ia_2v(t,x)^2\xi_1,   
\end{align*}
\begin{align*}
  B_1(t)
& =
  -
  a_0
  \Bigl\{
  (R_1\p_1\lvert{u}\rvert^2)
  +
  v[r_0(D),\bar{v}\p_1]
  +
  vR_1(\p_1\bar{v})+v\tilde{R}_1\bar{v}  
  \Bigr\}
\\
& -
  a_1\bar{u}\p_1u
  -
  a_2(u+v)\p_1\bar{u}
  -
  a_3(u+v)\bar{u},
\\
  B_2(t)
& =
  -
  a_0v
  \Bigl(
  [r_0(D),u\p_1]
  +
  R_1(\p_1u)+\tilde{R}_1u
  \Bigr)
  -
  a_1v\p_1u
  -
  a_3v^3. 
\end{align*}
By Lemma~\ref{theorem:syswp}, 
the initial value problem for the system of 
\eqref{equation:difference} and its complex conjugate 
is $L^2$-well-posed. 
Thus, ${}^t[w(t),\overline{w(t)}]=0$. 
\par
Lastly, we recover the continuity in the time variable. 
Let 
$u{\in}L^\infty(0,T;H^s(\mathbb{R}^2))$ 
be a unique solution to 
\eqref{equation:pde}-\eqref{equation:data}. 
Recall \eqref{equation:ryoko}. 
Let $\alpha$ be a multi-index satisfying 
$\lvert\alpha\rvert=[s]-1$. 
Set $\theta=s-[s]$ 
and 
$u_\alpha=\langle{D}\rangle^\theta\p^\alpha{u}$ 
for short. 
It suffices to show 
$u_\alpha{\in}C([0,T];H^1(\mathbb{R}^2))$. 
Applying $\langle{D}\rangle^\theta\p^\alpha$ to 
\eqref{equation:pde}, 
we have 
\begin{equation}
(\p_t+p_0(\p)+ip_1(\p)+a_1(t,x,D))u_\alpha 
+
a_2(t,x,D)\overline{u_\alpha}
=
f_\alpha,
\label{equation:minus}
\end{equation}
\begin{align*}
  a_1(t,x,\xi)
& =
  -
  i\lvert{u(t,x)}\rvert^2
  (a_0r_0(\xi)+a_1)\xi_1
  +
  ip_2(\xi),
\\
  a_2(t,x,\xi)
& =
  -
  iu(t,x)^2
  (a_0r_0(\xi)+a_2)\xi_1,
\end{align*}
\begin{align*}
  f_\alpha
& =
  \biggl\{
  \sum_{j=0}^3
  a_j
  \langle{D}\rangle^\theta\p^\alpha
  f_j(u)
\\
& \qquad
  -
  \lvert{u}\rvert^2(a_0R_1+a_1)\p_1
  \langle{D}\rangle^\theta\p^\alpha{u}
  -
  u^2(a_0R_1+a_1)\p_1
  \langle{D}\rangle^\theta\p^\alpha\bar{u}
  \biggr\}
\\
& +
  a_0\lvert{u}\rvert^2\tilde{R}_1\p_1
  \langle{D}\rangle^\theta\p^\alpha{u}
  +
  a_0u^2\tilde{R}_1\p_1
  \langle{D}\rangle^\theta\p^\alpha\bar{u}
\\
& -
  [\langle{D}\rangle^\theta,\lvert{u}\rvert^2(a_0R_1+a_1)\p_1]\p^\alpha{u}
\\
& -
  [\langle{D}\rangle^\theta,u^2(a_0R_1+a_1)\p_1]\p^\alpha\bar{u}
  -
  ip_3\langle{D}\rangle^\theta\p^\alpha{u}.
\end{align*}
It is easy to see that  
$u_\alpha(0){\in}H^1(\mathbb{R}^2)$ 
and 
$f_\alpha{\in}L^1(0,T;H^1(\mathbb{R}^2))$ 
in the same way as $f_{\ep,\alpha}$. 
It follows that 
$u_\alpha{\in}C([0,T];H^1(\mathbb{R}^2))$ 
since the initial value problem for the system of 
\eqref{equation:minus} and its complex conjugate 
is $H^1$-well-posed.   
%
%
\bibliographystyle{amsplain}

\end{document}